\documentclass[a4paper,10pt]{article}
\usepackage[cp1251]{inputenc}
\usepackage[english]{babel}

\usepackage{amssymb}\usepackage{amsmath}
\usepackage{tikz}
\usetikzlibrary{arrows,snakes,backgrounds}
\usepackage{graphicx}
\def\bl{\rule[-1mm]{2.4mm}{2.4mm}}

\def\be{\begin{equation}}
\def\ee{\end{equation}}

\def\ord{{\rm ord}}
\def\bl{\rule[-1mm]{2.4mm}{2.4mm}}
\def\vert{{\rule{.5mm}{2.mm}}}
\def\hor{{\rule{2.mm}{.5mm}}}

\newtheorem{thrm}{\bf Theorem}

\newtheorem{dfn}{\bf Definition}
\newtheorem{rmk}{\bf Remark}

\numberwithin{rmk}{section}
\numberwithin{lmm}{section}
\numberwithin{xmpl}{section}

\begin{document}
\title {Coordinate spaces of graphs: approaching interior faces}
\author{\copyright 2020 ~~~~A.B.Bogatyr\"ev
\thanks{Supported by RSF grant 16-11-10349}}
\date{}
\maketitle
{\bf Keywords:}{ Moduli space, real algebraic curve, Abelian integral, graphs, 
foliation of quadratic differential}\\
{\bf MSC 2010:} Primary 30F30, Secondary 32G15, 14H15\\
{\bf Abstract:} {\small
We consider the cell decomposition of the moduli space of real genus two curves with a marked point on the only real oval.
The cells are enumerated by certain graphs with their weights describing the complex structure on a curve.
We show that collapse of an edge of the graph results in a root like singularity of the natural mapping from the graph weights to 
the moduli space of curves.}

\section{Introduction}
Riemann surfaces do not necessarily appear as  algebraic curves. 
For some applications other representations are much more convenient. For instance, we can glue a surface
of some standard pieces of complex plane: half planes, triangles, rectangles, (half-) stripes, etc. The scars remaining after this surgery
make up a graph embedded to the surface. The idea to represent complex structures and even more fine objects (like abelian differentials,  quadratic differentials \cite{Str}, including Jenkins-Strebel ones, branched projective structures \cite{VF} etc.) on surfaces by embedded weighted graphs is not new. Possibly, Felix Klein \cite{Kl1,Kl2} was the initiator of the tradition.  The most prominent examples of this kind 
are \emph{Dessins d'Enfants} introduced by A.Grothendieck \cite{Gr,LZ}. Similar constructions were used by M.Bertola in his work on Boutroux curves \cite{Bert}. \emph{Ribbon graphs} once introduced by M.Kontsevich for the proof of Witten conjecture, today make up a separate flourishing industry inside mathematical physics, see e.g. \cite{Kon,Pen,Fock}.  The study of the period mapping lying in the core of the \emph{Chebyshev Ansatz} for the solutions of uniform rational approximation problems \cite{B02, Bbook} also grounds on a pictorial technique \cite{B03, B19}. Other closely related topics include \emph{Flat surfaces} \cite{Z}, \emph{Triangulated surfaces} \cite{VS};  see also \cite{LK, Sol} for other examples. Pictorial technique proved to be  extremely useful in the study of meso-scale geometry (as opposed to the differential and the global ones) of various moduli spaces and the related structures.

Graphs describing conformal  structures on a surface may be subjected to 
changes of their combinatorial characteristics like the flip transformation for the ribbon graph. 
It is intuitively clear that the complex structure of the underlying surface should have a continuous limit 
when we contract edges of the graph. However the detailed mathematical analysis of what is going on with the moduli
under this transformation has not been yet studied to the best of our knowledge. The response of moduli of a curve to the 
variations of weights of the graph strictly inside the space of admissible weights has been studied long ago. This dependence is 
real analytic, see e.g. \cite{Bbook,B03}. The same problem near the boundary of the admissible weights space is much more difficult since we have to compare the moduli of surfaces glued by different rules: this exactly reflects the change in combinatorial structure of the graph.

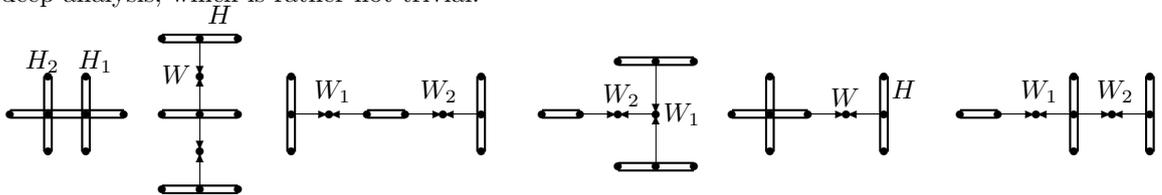
\begin{figure}[h!]
\begin{picture}(170,25)

\put(-5,15){
\begin{picture}(60,5)
\thicklines
\put(5,.5){\line(1,0){15}}
\put(5,-.5){\line(1,0){15}}
\put(10.5,5){\line(0,-1){10}}
\put(9.5,5){\line(0,-1){10}}
\put(15.5,5){\line(0,-1){10}}
\put(14.5,5){\line(0,-1){10}}
\thinlines

\multiput(5,0)(5,0){4}{\circle*{1}}
\put(10,5){\circle*{}}
\put(10,-5){\circle*{1}}
\put(15,5){\circle*{1}}
\put(15,-5){\circle*{1}}

\put(7,6){$H_2$}
\put(14,6){$H_1$}
\end{picture}}

\put(0,15){
\begin{picture}(60,10)
\thicklines
\put(20,.5){\line(1,0){10}}
\put(20,-.5){\line(1,0){10}}
\put(20,10.5){\line(1,0){10}}
\put(20,9.5){\line(1,0){10}}
\put(20,-10.5){\line(1,0){10}}
\put(20,-9.5){\line(1,0){10}}

\thinlines

\multiput(20,0)(5,0){3}{\circle*{1}}
\multiput(20,10)(5,0){3}{\circle*{1}}
\multiput(20,-10)(5,0){3}{\circle*{1}}
\put(25,5){\circle*{1}}
\put(25,-5){\circle*{1}}

\put(26,12){$H$}
\put(20,4){$W$}

\put(25,0){\vector(0,1){5}}
\put(25,10){\vector(0,-1){5}}
\put(25,0){\vector(0,-1){5}}
\put(25,-10){\vector(0,1){5}}
\end{picture}}

\put(17,15){
\begin{picture}(60,5)
\thicklines
\put(20.5,5){\line(0,-1){10}}
\put(19.5,5){\line(0,-1){10}}
\put(30,.5){\line(1,0){5}}
\put(30,-.5){\line(1,0){5}}
\put(45.5,5){\line(0,-1){10}}
\put(44.5,5){\line(0,-1){10}}
\thinlines

\multiput(20,0)(5,0){6}{\circle*{1}}
\put(20,5){\circle*{1}}
\put(20,-5){\circle*{1}}
\put(45,5){\circle*{1}}
\put(45,-5){\circle*{1}}

\put(20,0){\vector(1,0){5}}
\put(30,0){\vector(-1,0){5}}
\put(35,0){\vector(1,0){5}}
\put(45,0){\vector(-1,0){5}}

\put(23,2){$W_1$}
\put(37,2){$W_2$}

\end{picture}}

\put(50,15){
\begin{picture}(60,10)
\thicklines
\put(20,.5){\line(1,0){5}}
\put(20,-.5){\line(1,0){5}}
\put(30,7.5){\line(1,0){10}}
\put(30,6.5){\line(1,0){10}}
\put(30,-7.5){\line(1,0){10}}
\put(30,-6.5){\line(1,0){10}}

\thinlines

\multiput(20,0)(5,0){4}{\circle*{1}}
\multiput(30,7)(5,0){3}{\circle*{1}}
\multiput(30,-7)(5,0){3}{\circle*{1}}

\put(25,0){\vector(1,0){5}}
\put(35,0){\vector(-1,0){5}}

\put(35,7){\vector(0,-1){7}}
\put(35,-7){\vector(0,1){7}}

\put(36,-1){$W_1$}
\put(28,1.5){$W_2$}

\end{picture}}

\put(75,15){
\begin{picture}(60,5)

\thicklines
\put(25.5,5){\line(0,-1){10}}
\put(24.5,5){\line(0,-1){10}}

\put(40.5,5){\line(0,-1){10}}
\put(39.5,5){\line(0,-1){10}}
\put(20,.5){\line(1,0){10}}
\put(20,-.5){\line(1,0){10}}
\thinlines

\multiput(20,0)(5,0){5}{\circle*{1}}
\put(25,5){\circle*{1}}
\put(25,-5){\circle*{1}}
\put(40,5){\circle*{1}}
\put(40,-5){\circle*{1}}

\put(30,0){\vector(1,0){5}}
\put(40,0){\vector(-1,0){5}}

\put(41,2){$H$}
\put(33,1){$W$}

\end{picture}}

\put(105,15){
\begin{picture}(60,5)
\thicklines
\put(35.5,5){\line(0,-1){10}}
\put(34.5,5){\line(0,-1){10}}
\put(45.5,5){\line(0,-1){10}}
\put(44.5,5){\line(0,-1){10}}
\put(20,.5){\line(1,0){5}}
\put(20,-.5){\line(1,0){5}}
\thinlines

\multiput(20,0)(5,0){6}{\circle*{1}}
\put(35,5){\circle*{1}}
\put(35,-5){\circle*{1}}
\put(45,5){\circle*{1}}
\put(45,-5){\circle*{1}}

\put(25,0){\vector(1,0){5}}
\put(35,0){\vector(-1,0){5}}
\put(35,0){\vector(1,0){5}}
\put(45,0){\vector(-1,0){5}}

\put(28,2){$W_1$}
\put(38,2){$W_2$}

\end{picture}}
\end{picture}
\caption{\small  Graphs which together with their central symmetric graphs encode all full dimensional cells of the moduli 
space ${\cal H}_2^1$.  The  weights of vertexes and vertical edges we denote as $W$ and $H$ respectively.}
\label{H21}
\end{figure}

In this paper we consider a very concrete example: the moduli space ${\cal H}_2^1$ of genus two real curves with a unique real oval and a 
marked point on it disjoint from branch points. This space may be decomposed into nine full dimensional cells labeled by special trees shown in Fig \ref{H21}. 
We consider transition through a wall separating two neighbouring cells and show that natural moduli of curves 
given by the positions of branch points behave continuously but not smoothly with respect to the natural coordinates inside each cell:
root type singularity arises at the wall. In particular, the main result gives the justification of adjacency relations for the cell decompositions in the combinatorial  moduli space theory, which is often done (especially by physicists) at the intuitive level without deep analysis, which is rather not trivial.

The principal investigation method we use is the quasi-conformal \cite{Ahl} deformation of an abelian integral 
which is a simplification of the deformation technique proposed in \cite{Bbook}, Chap. 5. The quasi-conformal approach is
a traditional and a universal tool for the study of various  coordinate systems in moduli spaces, including the vicinity of the singularity 
of the coordinate change. The choice of the moduli space ${\cal H}_2^1$ as an object of study is to a certain extent accidental: this space is neither elementary, nor is it too involved, and the contraction of two types of edges can be illustrated at a time.
This space is also used for the analysis of certain applied problems \cite{B04, B05, Bbook}. The principal elements of our analysis can be extended to much more sophisticated cases.

{\bf Acknowledgements}. This research had been supported by Russian Scientific Foundation (Grant 16-11-10349)
and INM RAS division of Moscow Center for Fundamental and Applied Math (Agreement 075-15-2019-1624/2). The author thanks all the 
participants of A.Gonchar seminar on \emph{Complex analysis} (Steklov Math Inst.) and G.Shabat seminar '\emph{Graphs on surfaces and Curves over number fields}' (Moscow Lomonosov University) for the discussions related to this paper and for their constructive criticism. My special gratitude is to Prof. Hartmuth Monien for drawing attention to the works \cite{Kl1, Kl2} of F. Klein.

\section{Preliminary settings and main result}
In this section we remind the definition of the moduli space of real genus two curves with exactly one \emph{oriented} oval and a marked point on it; give the description of the curves by weighted graphs and formulate the main result of the paper.

\subsection{Moduli space}
There are two classes of smooth genus two real curves with a unique real oval: the oval contains either two branch (Weierstrass) points or no one of them. Let ${\cal H}_2^1$ be the moduli space of the curves $M$ from the first class, which are additionally equipped with a marked point $``\infty``$ (in what follows will be used without quotes) on the real oval. We require this marked point being not fixed by the hyperelliptic involution $J$ acting on each curve. This moduli space is used e.g. for the analysis of the problem about the so called optimal stability polynomials \cite{B04, B05, Bbook}, including the damped ones.

Any element of ${\cal H}_2^1$ admits the normalized affine model:  
\be
\label{M}
M=M({\sf E}):= \{(x,w)\in\mathbb{C}^2:\quad w^2=(x^2-1)\prod_{s=1}^2(x-e_s)(x-\bar{e}_s) \}, 
\ee
with branching points $e_1\neq e_2$ from the open upper half plane $\mathbb{H}$. Its branching set ${\sf E}=\{\pm1, e_s,\bar{e}_s\}_{s=1,2}$ 
has mirror symmetry $\sf E=\bar{E}$. The hyperelliptic and anticonformal involutions of the curve we define as
$J(x,w):=(x,-w)$ and $\bar{J}(x,w):=(\bar{x},\bar{w})$ respectively.
The marked point $\infty $ on the real oval is the point corresponding to $(x,w)=(+\infty,+\infty)$ in the natural two-point compactification of the affine curve (\ref{M}).  

The space ${\cal H}_2^1$ of such curves is parametrized by the positions of their branch points $e_1, ~e_2$ in the open upper half plane $\mathbb{H}$
which are indistinguishable and cannot coincide. Therefore we have a model $(\mathbb{H}^2\setminus\{diagonal\})/permutation$ for the moduli space
which shows that $\dim{\cal H}_2^1=4$ and  $\pi_1({\cal H}_2^1)=Br_2=\mathbb{Z}$.

\subsection{Distinguished differential}
On each curve $M$ from the moduli space there is a unique third kind abelian differential $d\eta_M$
with just two simple poles: the marked point $\infty$ and its involution $J\infty$, with residues $-1$ and $+1$ respectively and  purely imaginary periods \cite{Bbook}, \S 2.1.3. For the algebraic model (\ref{M}) of the curve $M$ the differential takes the form:
\be
\label{deta}
d\eta_M=(x^2+\dots)w^{-1}dx,
\ee
with dots standing for a linear polynomial. One can check that normalization conditions of the differential imply that the latter is \emph{real}, that is $\bar{J}d\eta_M=\overline{d\eta_M}$. In other words, the polynomial in (\ref{deta})
has real coefficients. An important consequence of this fact is this \cite{Bbook, B02}: the periods of this differential 
along even 1-cycles $C:=\bar{J}C$ vanish since they should be real and imaginary at the same time.

It is in terms of this differential that the solutions of various problems of uniform polynomial approximation may be represented
\cite{Bbook}:
$$
P_n(x)=\cos(ni\int_{(1,0)}^{(x,w)}d\eta_M),
$$
under the additional requirement $\int_{H_1(M,\mathbb{Z})}d\eta_M\subset\frac{2\pi i}n\mathbb{Z}$
which means that the periods of abelian integral are commensurable with the periods of $cos(\cdot)$ and therefore 
guarantees that the left hand side of the equality is a polynomial. Classical formulas for Chebyshev and Zolotarev 
polynomials are just particular cases of this representation.

\subsection{Global width function}
Suppose $M({\sf E})\in{\cal H}_2^1$ and $d\eta_M$ is the 3rd kind differential  associated with the curve $M$ as above.
One immediately checks that the normalization conditions of $d\eta_M$ imply that the \emph{width function} 
\be
\label{W}
W(x):=|Re\int_{(1,0)}^{(x,w)}d\eta_M|,
\qquad x\in\mathbb{C},
\ee
obeys the following properties:
\begin{itemize}
\item $W$ is single valued on the plane, 
\item $W$ is harmonic outside its zero set $\Gamma_\vert:=\{x\in\mathbb{C}: W(x)=0\}$, 
\item $W$ has a logarithmic pole at infinity,
\item $W$ vanishes at each branch point $e\in\sf E$\\
\end{itemize}

We only comment on the last property. Since $d\eta_M$ is odd with respect to the hyperelliptic involution of $M$,
$W(e_s)$ is equal to one half of the absolute value of the real part of some period of the differential.
Normalization implies that all its periods are purely imaginary.

The level lines of function  $W$ make up a vertical foliation of quadratic differential 
$(d\eta)^2$, whereas its steepest descent lines are horizontal trajectories  of the said differential.

\subsection{Construction of the graph $\Gamma(M)$.}
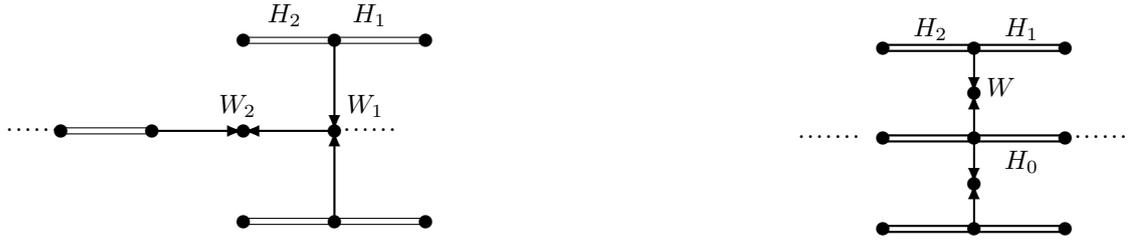
\begin{figure}
\begin{picture}(80,40)
\unitlength=.8mm
\thinlines
\put(10,20.5){\line(1,0){15}}
\put(10,19.5){\line(1,0){15}}
\multiput(40,35.5)(0,-1){2}{\line(1,0){30}}
\multiput(40,5.5)(0,-1){2}{\line(1,0){30}}
\thicklines
\put(25,20){\vector(1,0){15}}
\put(55,20){\vector(-1,0){15}}
\put(55,35){\vector(0,-1){15}}
\put(55,5){\vector(0,1){15}}

\multiput(10,20)(15,0){4}{\circle*{2}}
\multiput(40,35)(15,0){3}{\circle*{2}}
\multiput(40,5)(15,0){3}{\circle*{2}}

\multiput(10,20)(-1.5,0){7}{.}
\multiput(55,20)(1.5,0){7}{.}
\put(44,38){$H_2$}
\put(58,38){$H_1$}
\put(36,23){$W_2$}
\put(57,23){$W_1$}
\end{picture}
\begin{picture}(60,20)(-35,-15)
\unitlength=.8mm
\thicklines
\put(0,.5){\line(1,0){30}}
\put(0,-.5){\line(1,0){30}}
\put(0,15.5){\line(1,0){30}}
\put(0,14.5){\line(1,0){30}}
\put(0,-15.5){\line(1,0){30}}
\put(0,-14.5){\line(1,0){30}}


\multiput(0,0)(15,0){3}{\circle*{2}}
\multiput(0,15)(15,0){3}{\circle*{2}}
\multiput(0,-15)(15,0){3}{\circle*{2}}
\put(15,7.5){\circle*{2}}
\put(15,-7.5){\circle*{2}}

\put(20,17){$H_1$}
\put(20,-5){$H_0$}
\put(5,17){$H_2$}
\put(17,7){$W$}

\put(15,0){\vector(0,1){7.5}}
\put(15,15){\vector(0,-1){7.5}}
\put(15,0){\vector(0,-1){7.5}}
\put(15,-15){\vector(0,1){7.5}}

\multiput(-5,0)(-1.5,0){7}{.}
\multiput(30,0)(1.5,0){7}{.}
\end{picture}

\caption{\small Stable graphs $\Gamma_-$ (left) and $\Gamma_+$ (right) for curves of the moduli space ${\cal H}_2^1$. 
Double lines/arrows are vertical/horizontal edges; dotted line is the real (mirror symmetry) axis.
Independent weights $H$ of vertical edges and the values of the width function $W$ at the  vertexes of horizontal subgraph are shown.}
\label{StableG}
\end{figure}

To any curve $M$ from our moduli space we associate a weighted planar graph $\Gamma=\Gamma(M)$ composed of the finite number of segments of vertical and horizontal foliations \cite{Str} of the quadratic differential $(d\eta_M)^2$ descended to the Riemann sphere. The graph $\Gamma(M)$ is a union of the 'vertical' subgraph $\Gamma_\vert$ and the 'horizontal' subgraph $\Gamma_\hor$ -- see  Fig.\ref{StableG}  for examples of admissible graphs. 

\begin{dfn}
\begin{itemize} 
\item  Vertical edges are arcs of the  zero set of $W(x)$; they are segments of the vertical foliation $d\eta_M^2<0$ and are not oriented.  
\item Horizontal edges are all segments of the horizontal foliation $(d\eta_M)^2>0$ (or steepest descent lines for $W(x)$) 
connecting saddle points of function $W$ to other such points or -- as a rule -- to the zero set of $W$. Horizontal edges are oriented with respect to the growth of $W(x)$. 
\item Each edge, no matter what type is it, is equipped with its length in the metric $ds=|d\eta_M|$ of quadratic differential.
\item  Vertexes of the graph $\Gamma$ comprise all finite points of the divisor of the quadratic differential $({d\eta_M})^2$
considered on the plane as well as points in $\Gamma_\vert\cap\Gamma_\hor$ -- projections of the saddle points of $W$ to its zero set 
along the horizontal leaves. 
\end{itemize}
\end{dfn}

\begin{rmk} Instead of assigning lengths to the horizontal edges, it is more convenient to keep the values of the width function
$W(x)$ at all vertexes of the graph: the length of the oriented edge thus is the increment of the width function along it.
\end{rmk}

From the local behaviour of trajectories one immediately checks that the multiplicity of a vertex $V$ of the graph in the divisor of the quadratic differential $(d\eta_M)^2$ equals to the combinatorial value
$$
\ord (V):= d_\vert(V)+2d_{in}(V)-2,
$$
where $d_\vert$ is the degree of the vertex with respect to the vertical edges and $d_{in}$ is the number of incoming horizontal edges.
The branch points of the curve correspond to the vertexes $V$ of the subgraph $\Gamma_\vert\subset\Gamma$ with the odd value $\ord(V)$
-- see e.g. \ref{StableG}.  

A weighted planar graph may be associated to a (real) hyperelliptic curve with a marked point on it (its oval). 
The admissible graphs may be described in an axiomatic way: there are five restrictions on combinatorics and weights of the graphs
explicitly listed in \cite{B03, Bbook} which totally characterize them. The axioms may be actualized by a 
combinatorial algorithm listing all admissible graphs.

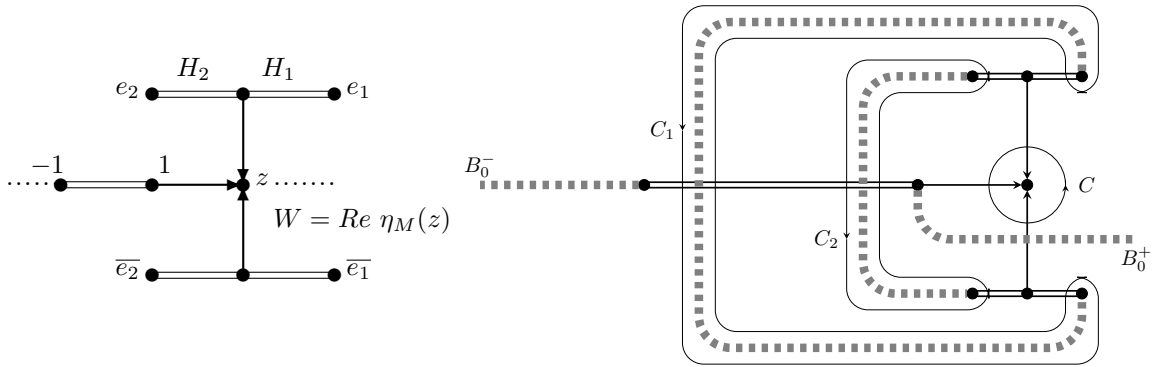
\begin{figure}
\begin{picture}(60,33)(0,-15)
\unitlength=.8mm
\thinlines
\put(10,20.5){\line(1,0){15}}
\put(10,19.5){\line(1,0){15}}
\multiput(25,35.5)(0,-1){2}{\line(1,0){30}}
\multiput(25,5.5)(0,-1){2}{\line(1,0){30}}
\thicklines
\put(25,20){\vector(1,0){15}}
\put(40,35){\vector(0,-1){15}}
\put(40,5){\vector(0,1){15}}

\multiput(10,20)(15,0){2}{\circle*{2}}
\put(40,20){\circle*{2}}
\multiput(25,35)(15,0){3}{\circle*{2}}
\multiput(25,5)(15,0){3}{\circle*{2}}

\multiput(10,20)(-1.5,0){7}{.}
\multiput(45,20)(1.5,0){7}{.}
\put(29,38){$H_2$}
\put(43,38){$H_1$}
\put(45,13){$W=Re~\eta_M(z)$}
\put(42,20){$z$}
\put(19,35){$e_2$}
\put(57,35){$e_1$}
\put(19,5){$\overline{e_2}$}
\put(57,5){$\overline{e_1}$}

\put(5,22){$-1$}
\put(26,22){$1$}
\end{picture}
\hfill
\includegraphics[trim = 4.8cm 10cm 5cm 11cm, clip, scale=.8]{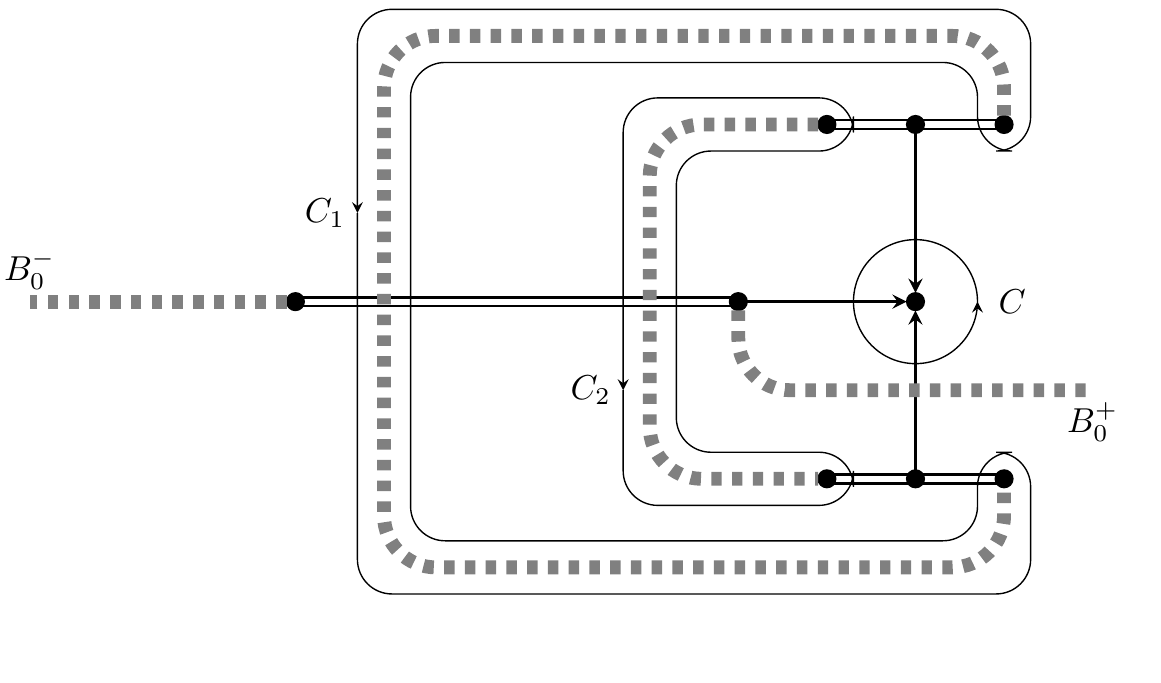}

\caption{\small Left:
Graph $\Gamma_{0}$  has codimension one coordinate space lying between full dimensional cells corresponding to graphs $\Gamma_+$ and $\Gamma_-$.
Right: Slit $B_s$, $s=1,2$, joining branch point $e_s$ to $\overline{e_s}$ (dashed thick gray line); contours $C_s$ encompass the slits $B_s$.
}
\label{UnStableG}
\end{figure} 
 
\subsection{Coordinate space of a graph}
It turns out that for the moduli space ${\cal H}_2^1$ there are only nine ``stable'' topological types of the graphs $\Gamma$ which do not change their combinatorial structure under arbitrary small perturbation of a curve $M({\sf E})$, they are listed e.g. in \cite{B19}. 
Two of such graphs $\Gamma_+$ and $\Gamma_-$  are shown in the left and the right panels of Fig. \ref{StableG}. The admissible independent weights sweep out four-dimensional space ${\cal A}[\Gamma]$ called coordinate space of a graph admitting  explicit description 
\cite{B03, Bbook}:
$${\cal A}[\Gamma_-]=\{(H_1,~H_2,~W_1,~W_2)\in\mathbb{R}_+^4: \quad 2(H_1+H_2)<\pi;\quad W_1<W_2\},$$ 
$${\cal A}[\Gamma_+]=\{(H_0,~H_1,~H_2, ~W)\in\mathbb{R}_+^4:\quad H_0+2(H_1+H_2)<\pi\},\quad\mathbb{R_+}:=(0,\infty).$$
Each coordinate space is the interior of the product of a symplex swept out by variables $H$ by a cone swept by variables $W$.

Along with topologically stable types of graphs $\Gamma$, there are lots of intermediate unstable
ones like graph $\Gamma_0$  shown  of the Fig. \ref{UnStableG}. Its coordinate space has codimension one in the moduli space ${\cal H}_2^1$: 
$${\cal A}[\Gamma_0]=\{(H_1,~H_2, ~W)\in\mathbb{R}_+^3:\quad 2(H_1+H_2)<\pi\},$$
and serves as an interface between coordinate spaces of $\Gamma_\pm$.

Each of the graphs $\Gamma_-, \Gamma_0, \Gamma_+$ with weights from respective coordinate space may be realized as a graph of a unique curve 
$M$ from the moduli space ${\cal H}_2^1$ \cite{B03, Bbook}. This Riemann surface may be glued from a finite number of stripes
in a way determined by combinatorics and weights of the graph. The detailed instruction for this assembly is given in \cite{B03, Bbook, B19}.   
Unfortunately, this approach to the reconstruction of the curve cannot be called efficient and it does not give any quantitative 
characteristics of the embedding of a coordinate space to the moduli space. We will use more flexible and constructive methods related to quasi-conformal mappings in what follows.

This author claimed \cite{B19} that the embedding of a positive codimension coordinate space like
${\cal A}[\Gamma_0]$ coincides with the continuation of the embedding of suitable full dimensional polyhedra ${\cal A}[\Gamma]$
to their faces. Here we show the validity of this statement in a particular case. 
The embedding of a coordinate space to the moduli space \cite{Bbook, B03} is 
given by a pair of complex valued real analytic functions $e_s(H,W),\quad s=1,2,$  
defined inside a polyhedron ${\cal A}[\Gamma]$ and describing the dependence of the branching divisor of the curve of the weights of the graph 
$\Gamma$.  In this paper we study the behaviuor of those functions near the 
boundary $\{W_1=W_2\}$ of the coordinate space ${\cal A}[\Gamma_-]$ and the boundary $\{H_0=0\}$ of the polyhedron ${\cal A}[\Gamma_+]$.

\subsection{Main theorem}
Consider an arbitrary point $A^0:=(H_1,~H_2,~W)$ of codimension one coordinate space ${\cal A}[\Gamma_0]$ and its small 
displacement $\delta A:=(~\delta H_1, ~\delta H_2, ~\delta W)$ within this space. The ``transversal''  displacement to the neighbouring 
full dimensional spaces ${\cal A}[\Gamma_\pm]$ will be described by a (small) positive variable $h$ which together with the tangential shift
$\delta A$ defines two points:
\be
\label{disp-}
{\cal A}[\Gamma_-]\ni\quad A^-=(H_1^-,H_2^-,W_1^-,W_2^-):=(H_1+\delta H_1, H_2+\delta H_2, W+\delta W-2h^3, W+\delta W+2h^3),
\ee
\be
\label{disp+}
{\cal A}[\Gamma_+]\ni\qquad A^+=(H_0^+,H_1^+,H_2^+,W^+):=(2h^3, H_1+\delta H_1 -2h^3, H_2+\delta H_2+2h^3, W+\delta W);
\ee
To each of points $A^\pm$ of coordinate spaces there corresponds a normalized branching divisor ${\sf E}\in{\cal H}_2^1$ containing two points $e\neq\pm1$ in the upper half plane. 

\begin{thrm}\label{mainThrm}
The displacement of the branch point $e$ caused by the tangential displacement $\delta A$ and the transversal displacement
$h$ admits the following asymptotic expansion:

\be\label{asymp}
2\pi i~(e(A^\pm)-e(A^0)) = 
\left\{
i\delta H_1\int_{C_1} - i\delta H_2\int_{C_2}  +
\delta W\int_C  \right\}
d\eta^e \pm 3h^2 \int_C y(x)d\eta^e + O((|\delta A|+h^2)^2),
\ee
where the meromorphic differential $d\eta^e$ on the unperturbed curve $M=M(A^0)$ is defined by the equality
$d\eta^ed\eta_M:=\displaystyle{\frac{e^2-1}{x^2-1}\frac{(dx)^2}{x-e}}$ with  finite area quadratic differential holomorphic in $\mathbb{CP}^1\setminus{\sf E}$
in the r.h.s.; $C,C_1$ and $C_2$ are even cycles on the curve $M$ shown in the right panel of Fig. \ref{UnStableG};
$y(x):=\alpha(x-z)+\dots$ is the real local coordinate on the curve in the vicinity of the double zero $z$ of the distinguished differential $d\eta_M$, defined by the equality $\eta_M(x):=W+y^3$; $|\delta A|$ is the euclidean length of $\delta A$.
\end{thrm}
 
\begin{rmk}
\label{mainTh}
(i) We observe that the embedding of the coordinate space to the moduli space contains root-like singularity 
$(H_0)^{2/3}$ for the space ${\cal A}[\Gamma_+]$ and $(W_2-W_1)^{2/3}$ for the space ${\cal A}[\Gamma_-]$ near appropriate 
boundary of the space and in particular it is not continuously  differentiable up to the boundary.

(ii) Two integrals around the pole $z$ in formula \eqref{asymp} may be calculated explicitly by residues:
\be
\begin{array}{l}
\int_C d\eta^e=2\pi i(\frac{\Omega'(z)}{3\alpha^3}-\frac{4\beta^4\Omega(z)}{9\alpha^6})\\
\int_C y(x)d\eta^e=2\pi i\frac{\Omega(z)}{3\alpha^2}
\end{array}
\ee
where $\Omega(z)=\frac1{z-e}\frac{e^2-1}{z^2-1}$ is the coefficient of the quadratic differential $d\eta_Md\eta^e$ and $\alpha>0, \beta$ are coefficients of the expansion
$\eta_M(x)=W+\alpha^3(x-z)^3+\beta^4(x-z)^4+\dots$.
\end{rmk}

\section{Main theorem proof}
Plan of the proof is as follows:  we consider the distinguished abelian integral $\eta_M(x):=
\int_{(1,0)}^{(x,w)}d\eta_M$ on the unperturbed curve $M$ corresponding to a chosen point $A^0$
from the coordinate space of the unstable graph $\Gamma^0$. To make it single valued on the complex plane we have to introduce three slits 
pairwise joining the branch points of the curve. Then, given sufficiently small tangential $\delta A$ and transversal $h$ displacements 
we explicitly construct (a) a tiny deformation $\eta^\pm(x)$ of the function $\eta_M$ and (b) a new global variable $\xi(x)$ in the complex plane, such that $\eta^\pm$ considered as a function of the new coordinate will be the distinguished abelian integral on a modified curve $M^\pm$ corresponding to a point $A^\pm$ in the coordinate space of the  stable graph $\Gamma_\pm$. Ahlfors formula for infinitesimal quasiconformal mapping will eventually give us the asymptotic formulas \eqref{asymp} for the displacements of the branch points $e\in{\sf E}$
of the curve. Now we proceed to the step by step realization of this plan.

\subsection{Abelian integral} The  integral $\eta:=\eta_M$ of the distinguished differential on the chosen curve $M$ corresponding to
the point $A^0$ of the codimension one coordinate space ${\cal A}[\Gamma_0]$ is locally single valued function of $x$ with the exception of 
infinity and the points of $\sf E$ where it has branching. We introduce three disjoint slits $B_j$, $j=0,1,2$ connecting 
the branch points pairwise, one of those passes through infinity and we designate two parts of it as 
$B_0^-$ and $B_0^+$  -- as in the  right Fig.  \ref{UnStableG}. The abelian integral admits a single-valued branch $\eta(x)$ in the remaining 
3-connected domain (pants) in the complex plane of variable $x$. Indeed, the homology basis of pants
(say the contours $C_1$ and $C_2$) are lifted to the \emph{even} cycles of the curve $M$, which by definition endure the action of the reflection $\bar{J}$. The integral of real differential $d\eta_M$ over even cycle is real (due to mirror symmetry) on the one hand and purely imaginary (due to normalization) on the other hand, hence it is zero.

Note that the sum of boundary values of the abelian integral $\eta$ is locally constant along each slit since the 
distinguished differential $d\eta_M$ is odd with respect to hyperelliptic involution $J$. The values of those constants are 
purely imaginary due to the normalization of $d\eta_M$. They may be easily calculated since we can reconstruct 
the values of the integral in all remarkable points related to the curve, given the weights of the graph $\Gamma^0$ (this calculation requires careful consideration of signs):
\be
\label{etaVal}
\begin{array}{l}
\eta(e_1)=\eta(\overline{e_1})=iH_1; \qquad
\eta(e_2)=\eta(\overline{e_2})=-iH_2;\\
\eta(1)=0;\qquad \eta(z)=W; \qquad \eta(-1)=i\pi. 
\end{array}
\ee

\subsection{Deformation of abelian integral}
We consider two smooth deformations of the abelian integral $\eta$ caused by the displacement $\delta A:= (\delta H_1,~ \delta H_2,~ \delta W)$
in the space ${\cal A}[\Gamma_0]$ and the transversal displacement $h>0$,  which correspond to the choice of the sign $\pm$
in the following formula:
\be
\label{deformation}
\eta^\pm(x):=\eta(x) -i\sum_{s=1,2}(-1)^s\delta H_s\rho_s(x)+(\delta W \pm 3h^2y(x))\rho(x).
\ee
Here $0\le\rho_s(x)\le1$ is a smooth real valued cut-off function equal to $1$ in a  vicinity of the arc $B_j$ 
and vanishing identically outside some larger vicinity of the same arc; $\rho(x)$ is a similar function equal to $1$ in the 
vicinity of the double zero $z$ of the distinguished differential $d\eta_M$. The supports of different
cut-off functions do not intersect. The real coordinate $y(x)=\alpha(x-z)+\dots$ in the vicinity of $z$ is specified in the formulation of theorem \ref{mainThrm}.

\subsection{New global coordinate}
Let $\xi(x)=\xi(x;~\delta A, h, \pm)$ be the  solution of the Beltrami equation
$\xi_{\bar{x}}=\mu(x)\xi_x$ with the coefficient $\mu(x):=\eta^\pm_{\bar{x}}/\eta^\pm_x$, homeomorphic in the whole plane
and pinning three points $x=\pm1,\infty$. The support of Beltrami coefficient $\mu$ 
consists of three annular domains encompassing the slits $B_1$, $B_2$ and the double pole $z$ of  differential $d\eta^e$. 

We suppose w.l.o.g. that the values of each cut-off function $\rho_*(x)$   coincide in complex conjugate points. This implies that 
Beltrami coefficient will be mirror symmetric: $\mu(\bar{x})=\bar{\mu}(x)$. Since the normalizing set $x=\pm1,\infty$ of the map
is real, the new variable $\xi(x)$ will be mirror symmetric (real) too. 

\subsection{Key observation}
We claim  that the perturbed function $\eta^\pm(x(\xi))$ considered as a function of the new global variable $\xi$ 
is the distinguished abelian integral for some disturbed curve $M^\pm$ with the branching divisor ${\sf E}^\pm:=\xi({\sf E};~\delta A,h,\pm)$
parametrically dependent on the displacements. 

First, we check that $\eta^\pm$ is a holomorphic function of variable $\xi$ outside the system of cuts
$\xi(B_s)$, $s=0,1,2$. Indeed, the inverse mapping $x(\xi)$ is also  quasi-conformal with Beltrami coefficient $\nu(\xi)$ satisfying the relation \cite{Ahl}: 
$$
\nu(\xi)x_\xi+\mu(x)\overline{x_\xi}=0,
$$
which is obtained by differentiating the identity $\xi(x(\xi))=\xi$. Now
$$
\eta^\pm_{\bar{\xi}}=\eta^\pm_x x_{\bar{\xi}}+\eta^\pm_{\bar{x}} \bar{x}_{\bar{\xi}}=
\eta^\pm_x(x_{\bar{\xi}}+\mu(x)\bar{x}_{\bar{\xi}})=\eta^\pm_x(x_{\bar{\xi}}-\nu(\xi)x_\xi)=0.
$$

Next step is to check that the boundary values of $\eta^\pm$ on the banks of the cuts sum up to a purely imaginary constant, individual for every cut.  This is clearly seen from the formula \eqref{deformation}: this constant equals to $-2i(H_s+\delta H_s)(-1)^s$ for the cut $B_s$, $s=1,2$; $0$ for $B_0^+$ and $i\pi$ for $B_0^-$. 
This observation implies $(d\eta^\pm)^2$ being a rational quadratic differential on the sphere of variable $\xi$ with simple poles at the points of $\xi({\sf E})$ and double pole at infinity. The residue  of $d\eta^\pm$  at infinity is the same as for the differential 
$d\eta_M$ of the unperturbed curve. Pure imaginarity of constants found above ensure real normalization of the new differential.

\subsection{Infinitesimal deformation}
We shall use Ahlfors formula for the infinitesimal quasi-conformal deformations to obtain the low order terms in deformation of the branching divisor $\sf E$. First we assess the Beltrami coefficient keeping the leading terms of the deformation  only:
$$
\eta^\pm_{\bar{x}}=-i\sum_{j=1,2}(-1)^j\delta H_j \rho_{j\bar{x}}(x) \quad+ (\delta W\pm 3h^2 y(x)) \rho_{\bar{x}}(x);
$$
$$
\eta^\pm_x=\eta_x-i\sum_{j=1,2}(-1)^j\delta H_j \rho_{jx}(x) \quad+ (\delta W\pm 3h^2 y(x)) \rho_x(x)\pm3\rho h^2 y_x(x);
$$
\be\label{mu}
\mu(x)=-i\sum_{j=1,2}(-1)^j\frac{\rho_{j\bar{x}}(x)}{\eta_x}\delta H_j + \frac{\rho_{\bar{x}}(x)}{\eta_x}(\delta W\pm3h^2y(x)) +
O(|\delta A|^2+ h^4+|\delta A|h^2)
\ee
with uniformly bounded residual term on the support of $\mu$.

For each branch point $e\in\sf E$ we have the expression of its displacement caused by the change of moduli in the unstable coordinate space 
and the transversal displacement to the stable space:
\be
\label{deltaE}
2\pi i~\delta e= \int_{Supp~\mu} \frac{e^2-1}{x^2-1}\frac{\mu(x)}{x-e} dx\wedge\bar{dx} +O(||\mu||^2_\infty).
\ee
Now we insert the approximate expression for $\mu(x)$ from \eqref{mu} to this formula and use integration by parts
which reduces the said formula to
$$
2\pi i~\delta e^\pm=\int_{\partial~ Supp~\mu}
i\sum_{j=1,2}(-1)^j\delta H_j\rho_j(x)d\eta^e-
(\delta W\pm3y(x)h^2)\rho(x)d\eta^e+O(\dots),
$$
where meromorphic differential $d\eta^e:=(e^2-1)\frac{(x-\bar{e})(x-e')(x-\bar{e'})dx}{(x-z)^2w}$, $e'\in{\sf E}\setminus\{\pm1,e\}$, was 
defined in the formulation of the theorem \ref{mainThrm} and
the order of magnitude of the residual term is the same as in \eqref{mu}. 
The value of each mollifier $\rho_*$ is equal to 1 on exactly one contours of six which bound the support of Beltrami coefficient
and vanishes on the rest five, therefore the displacement of a branch point takes the form \eqref{asymp}. 

\subsection{New weighted graph}
It is easy to draw the graph $\Gamma$ corresponding to the deformation \eqref{deformation} of the distinguished abelian integral in $x$-coordinate. The value of the function $W(x)=Re~\eta^\pm(x)$ is changed in the vicinity of the point $z$ only, therefore the vertical part $\Gamma_\vert$ of the graph remains intact. The horizontal part will be slightly modified since the deformation splits the double critical point $z$ into two simple ones.

To find zeros of the disturbed differential $d\eta^\pm$ we use  local coordinate $y$ in the vicinity of the double zero 
$z$ of the undisturbed differential $d\eta_M$. Putting $\rho=1$ in this vicinity we get the representation 
$\eta^\pm=W+\delta W+ y^3\pm 3h^2y$ for the integral. For small enough $h>0$ zeros $y=\pm h$ (for the deformation labeled "-")
or $y=\pm ih$ (for the deformation labeled "+") of $d\eta^\pm$ lie in the same vicinity of the double zero with coordinate $y(z)=0$. 
Drawing the isolines of the function $Im~\eta^\pm$ passing through the found critical points till their intersection with $\Gamma_\vert$
we see that the graph associated to the new curve $M^\pm$ is  $\Gamma_\pm$.

The critical values of the deformation of abelian integral together with the values at other remarkable points 
allow us to reconstruct the weights of the graph $\Gamma_\pm$.
Take for example the deformation "+".  For an arbitrary point $A^+=(H_0^+,H_1^+,H_2^+,W^+)\in{\cal A}[\Gamma_+]$ the critical values of the function $\eta_M^+$ defined outside the system of cuts homotopic to $B_s$, $s=0,1,2$, are equal to $W^+\pm iH_0^+$. The values of $\eta_M^+$
at the branch points $e_1,e_2$ are $i(H_0^++H_1^+)$ and $i(H_0^+-H_2^+)$ respectively (Cf.:\eqref{etaVal}). Calculations take into account the fact that branches of distinguished differential outside the graph $\Gamma_+$ and outside the system of cuts $B_s$ differ at most by sign. 
Comparison to the values explicitly taken from \eqref{deformation}, e.g. critical values 
$\eta^+(\pm ih)=W+\delta W\pm 2ih^3$, shows that the new curve $M^+$ has coordinates $A^+=(2h^3, H_1+\delta H_1 -2h^3, H_2+\delta H_2+2h^3, W+\delta W)$ in the space ${\cal A}[\Gamma_+]$. 

Similar calculations for the case of deformation "-" bring us to the point \eqref{disp-} of the coordinate space ${\cal A}[\Gamma_-]$.
~~\bl

\section{Conclusion}
We have studied the dependence of the branch points of a curve from the length of the vanishing 
edge of a graph describing the conformal structure.  To find the asymptotic  of the embedding of polyhedra ${\cal A}[\Gamma^\pm]$
near the face corresponding to the vanishing edge, we have adapted the technology of quasi-conformal 
deformation of abelian integrals elaborated in \cite{Bbook}, Chap.5. It turned out that the dependence has a '\emph{cuspidal}' singularity  
with respect to the vanishing transversal coordinate of the polyhedron and smooth with respect to all tangential coordinates.
All constants in the obtained asymptotic are given in explicit controllable form as periods of a certain abelian integral.
Higher order terms of the expansion \eqref{asymp} may also be explicitly calculated and will be presented elsewhere.
Similar expansions may be given in the vicinity of the ``exterior'' faces of the coordinate space
${\cal A}[\Gamma_\pm]$ corresponding to the merger of two branch points.

\vspace{5mm}
\parbox{9cm}
{\it
119991 Russia, Moscow GSP-1, ul. Gubkina 8,\\
Institute for Numerical Mathematics,\\
Russian Academy of Sciences\\[2mm]
{\tt ab.bogatyrev@gmail.com}}

\begin{thebibliography}{80}
\bibitem{Ahl} L.Ahlfors, Lectures on Quasiconformal mappings -- AMS, University lecture series, 2006.
\bibitem{Bert} Bertola M., Boutroux curves with external field: equilibrium measures without a
minimization problem, Anal. Math. Phys. 1 (2011), 167-211.
\bibitem{B02} Bogatyrev A.B., Effective approach to least deviation problems.//Sbornik: Math  193:12 (2002), 1749-1769.
\bibitem{B03} Bogatyrev A.B., Combinatorial description of a moduli space of curves and of Extremal Polynomials //Sbornik: Math 194:10 (2003), 1451--1457.  See also: Errata to Bogatyrev A.B. Combinatorial description of a moduli space of curves and of Extremal Polynomials //Sbornik: Math 194:12 (2003), 1899.
\bibitem{B04}  	A.B. Bogatyrev, Effective computation of optimal stability polynomials, Calcolo, 41:4 (2004), 247--256
\bibitem{B05} A.B. Bogatyrev, Effective solution of the problem of the optimal stability polynomial, Sb. Math., 196:7 (2005), 959--981 
\bibitem{Bbook} Bogatyrev A.B. Extremal Polynomials and Riemann Surfaces  -- MCCME, 2005 (in Russian), Springer translation, 2012. 
\bibitem{B19} Bogatyrev A.B., Combinatorial description of period mapping: topology of 2D fibers// Sbornik: Math, 210:11 (2019)  
\bibitem{Fock} L.O.Chekhov and V.V. Fock, Quantum Teichmueller spaces// Theor. Math. Phys., 120:3, 511--528 (1999)
\bibitem{VF} V.V. Fock,  Description of moduli space of projective structures via fat graphs, arXiv:hep-th/9312193 
\bibitem{Gr} Alexandre Grothendieck, Esquisse d'un Programme// Geometric Galois Actions. -- Cambridge: Cambridge University Press., pp. 7--48.
\bibitem{Kl1} Klein F., \"Uber die Ernidrigung der Modulargleichungen //Math. Ann., Bd 14, 1879.
\bibitem{Kl2} Klein F., \"Uber die Transformation elfter ordnung der elliptischen Funktionen// Math. Ann., Bd 15, 1879.
\bibitem{Kon} M. L. Kontsevich, Intersection theory on the moduli space of curves//Func. Anal. and Appl., 25:2 (1991),  50--59.
\bibitem{LZ} S.K. Lando, A.K.Zvonkine, Graphs on surfaces and their applications, Springer 2004. 
\bibitem{LK} S. K. Lando, I.M.Krichever and A.S.Skripchenko, Foliations on spaces of real-normalized differentials,  arXiv:2010.09358  
\bibitem{Pen} R.C.Penner,  The decorated Teichmueller space of Riemann surfaces//Comm. Math. Physics, 113:2, 299--339, (1987)
\bibitem{Sol} A.Yu.Solynin, Quadratic differentials and weighted graphs on compact surfaces// Analysis and Math Physics
(B.Gustafsson and A.Vasil'ev eds.), Birkhauser, 2009  
\bibitem{Str} K.Strebel, Quadratic differentials -- Springer 1984  
\bibitem{VS} V. A. Voevodsky, G. B. Shabat, Equilateral triangulations of Riemann surfaces, and curves over algebraic number fields, Dokl. Math., 39:1 (1989), 38--41 
\bibitem{Z} A.Zorich, Flat surfaces, Frontiers in number theory, physics, and geometry. I, 437-583, Springer, Berlin, 2006;
 arXiv:math/0609392

 



\end{thebibliography}
\end{document}